\documentclass[12pt]{article}
\usepackage{amsmath} 
\usepackage{amsthm} 
\usepackage{amssymb} 
\usepackage[dvipsnames]{xcolor} 
\usepackage{mathtools} 
\usepackage{mma}
\usepackage{tikz,tikz-cd}
\usepackage[cal=boondoxo,scr=euler]{mathalfa} 
\usepackage[backref=page,linktocpage]{hyperref} 
\usepackage{graphics,graphicx} 
\usepackage{cleveref}
\usepackage{enumerate} 
\usetikzlibrary{calc}

\DeclareMathAlphabet{\mathsf}{OT1}{\sfdefault}{m}{n} 
\usepackage[margin=1.20in]{geometry} 

\hypersetup{ 
	colorlinks = true,
	linkbordercolor = {white},
	linkcolor = {BrickRed},
	anchorcolor = {black},
	citecolor = {BrickRed},
	filecolor = {cyan},
	menucolor = {BrickRed},
	runcolor = {cyan},
	urlcolor = {black}
}

\newtheoremstyle{teoremas} 
{11pt}
{11pt}
{\itshape}
{}
{\bfseries}
{}
{.5em}
{}

\theoremstyle{teoremas} 
\newtheorem{theorem}{Theorem}[section] 
\newtheorem{corollary}[theorem]{Corollary} 
\newtheorem{lemma}[theorem]{Lemma} 
\newtheorem{proposition}[theorem]{Proposition} 

\newtheoremstyle{definition} 
{11pt}
{11pt}
{}
{}
{\bfseries}
{}
{.5em}
{}

\theoremstyle{definition} 
\newtheorem{conjecture}[theorem]{Conjecture} 
\newtheorem{remark}[theorem]{Remark} 

\crefname{theorem}{theorem}{theorems} 
\Crefname{theorem}{Theorem}{Theorems} 
\crefname{lemma}{lemma}{lemmas} 
\Crefname{lemma}{Lemma}{Lemmas} 
\crefname{proposition}{proposition}{propositions} 
\Crefname{proposition}{Proposition}{Propositions} 

\DeclareMathOperator{\rk}{rk}
\DeclareMathOperator{\cusp}{cusp}
\DeclareMathOperator{\Rel}{Rel}

\newcommand{\M}{\mathsf{M}}
\newcommand{\PP}{\mathscr{P}}

\newcommand{\U}{\mathsf{U}}

\newcommand{\intt}{\operatorname{int}}
\newcommand{\LL}{\mathsf{\Lambda}}

\renewcommand{\dim}{\operatorname{dim}}

\makeatletter
\def\dual#1{\expandafter\dual@aux#1\@nil}
\def\dual@aux#1/#2\@nil{\begin{tabular}{@{}c@{}}#1\\#2\end{tabular}}
\makeatother

\begin{document}

\begin{center}
{\large \bf The inverse $Z$-polynomial of a matroid}
\end{center}

\begin{center}
 Alice L.L. Gao$^{1}$, Xuan Ruan$^{2}$ and Matthew H.Y. Xie$^{3}$\\[6pt]

$^{1,2}$Shenzhen Research Institute of Northwestern Polytechnical University,\\
Sanhang Science $\&$ Technology Building, No. 45th, Gaoxin South 9th Road, Nanshan District, Shenzhen City, 518057, P.R. China

$^{1,2}$School of Mathematics and Statistics,\\
 Northwestern Polytechnical University, Xi'an, Shaanxi 710072, P.R. China

 $^{3}$College of Science, \\
 Tianjin University of Technology, Tianjin 300384, P. R. China\\[6pt]

 Email: $^{1}${\tt llgao@nwpu.edu.cn},
 $^{2}${\tt ruanxuan1211@163.com},
 $^{3}${\tt xie@email.tjut.edu.cn}
\end{center}

\noindent\textbf{Abstract.}
Motivated by the $Z$-polynomials of matroids, Ferroni, Matherne, Stevens, and Vecchi introduced the inverse $Z$-polynomial of a matroid. In this paper, we prove several fundamental properties of the inverse $Z$-polynomial, including non-negativity and multiplicativity, and show that it is a valuative invariant. We also provide explicit formulas for the inverse $Z$-polynomials of uniform matroids and a broader class of matroids, namely sparse paving matroids, which include uniform matroids as a special case. Furthermore, we establish the unimodality and log-concavity of these polynomials in the case of sparse paving matroids. Based on the properties of the $Z$-polynomial, we conjecture that the coefficients of the inverse $Z$-polynomial are unimodal and log-concave.

\noindent \emph{AMS Classification 2020:} 05B35, 52B40, 05A20


\noindent \emph{Corresponding Author:} Matthew H.Y. Xie, xie@email.tjut.edu.cn


\section{Introduction}

Since the introduction of Kazhdan-Lusztig polynomials of matroids by Elias, Proudfoot, and Wakefield \cite{elias2016kazhdan}, these polynomials have attracted 
much attention. For instance, see \cite{BV_2020EJC,ferroni2024braid,ferroni2023stressed,ferroni2022valuative,ferroni2022matroid,gao2021uniform,gao2024braid,gao2021inverse,Gedeon_2017ejc,Gedeon_proud_2017sem,proudfoot2018algebraic} and references therein. 
In \cite{PXY_2018ELC}, Proudfoot, Xu, and Young introduced the $Z$-polynomial $Z_{\M}(t)$ of a matroid $\M$, which is defined as the weighted sum of the Kazhdan-Lusztig polynomials for all possible contractions of $\M$. 
They also demonstrated that, for realizable matroids, the coefficients of $Z_{\M}(t)$ can be interpreted as the intersection cohomology Betti numbers of a projective variety. 
A precise description of this intersection cohomology for any matroid was later provided by Braden, Huh, Matherne, Proudfoot, and Wang \cite{BHM_ARX}.
Within the framework of Kazhdan-Lusztig-Stanley theory, Ferroni, Matherne, Stevens, and Vecchi \cite{FMSV_av_2024} introduced the inverse $Z$-polynomial of a matroid. 
Their work in \cite{FMSV_av_2024} primarily focuses on the Hilbert-Poincar\'e series of the Chow rings and augmented Chow rings of matroids and 
does not pay much attention to the inverse $Z$-polynomials, despite their many interesting properties that are worth further exploration.
The main objective of this paper is to compute the inverse $Z$-polynomials of uniform matroids and sparse paving matroids. In addition, we study the unimodality and log-concavity of the inverse $Z$-polynomials of sparse paving matroids.

Let us first recall some basic definitions related to the inverse $Z$-polynomials of matroids. 
For readers unfamiliar with matroids, 
we recommend Oxley’s book \cite{ox2011matroid} on matroid theory for background.
Throughout this paper, we adopt the notions 
and notation used in Ferroni, Matherne, Stevens, and Vecchi \cite{FMSV_av_2024}.
Let $\M = (E, \mathscr{B})$ denote a matroid, 
where $E$ is the ground set and $\mathscr{B}$ is the set of all bases of $\M$. The rank of a subset $A \subseteq E$ is defined by
$\rk(A) := \max_{B \in \mathscr{B}} |A \cap B|$, and the rank of $\M$ is given by $ \rk(\M) := \rk(E).
$
A flat of $\M$ is a subset $F \subseteq E$ such that either $F = E$ or the rank of $F$ strictly increases upon adding new elements. 
The set of all flats of $\M$, ordered by inclusion, forms a lattice, denoted by $\mathcal{L}(\M)$.
For any subset $A \subseteq E$, the restriction of $\M$ to $A$ is denoted $\M|_A$, and the contraction of $\M$ by $A$ is denoted $\M/A$. 
Let $P_{\M}(t)$ and $Q_{\M}(t)$ denote the Kazhdan-Lusztig and inverse Kazhdan-Lusztig polynomial of $\M$, respectively. 
Given a loopless matroid $\M$, the $Z$-polynomial of $\M$ is defined as
\begin{align}\label{z_defini}
Z_\M(t):\,=\sum_{F\in \mathcal{L}(\M)}
t^{\mathrm{rk}(\M|_F) } P_{\M/F}(t)
\end{align}
by Proudfoot, Xu, and Young \cite{PXY_2018ELC}.
Working in the incidence algebra framework, Ferroni, Matherne, Stevens, and Vecchi \cite[p.25]{FMSV_av_2024} introduced a related polynomial $\hat{Y}_\M(t)$, defined by 
\begin{align}\label{z-inverse-algebra}
\hat{Y}_\M(t):\,=\sum_{F\in \mathcal{L}(\M)}
 (-1)^{\mathrm{rk}(\M|_F)} Q_{\M|_F}(t) \cdot t^{\mathrm{rk}(\M/F) } \mu_{\M/F},
\end{align}
where $\mu_{\M}=\mu[\emptyset,E]$ is the 
M\"{o}bius invariant of $\M$.
We define the \emph{inverse $Z$-polynomial} of $\M$ by
\begin{align}\label{defi-z-mod}
Y_{\M}(t) := (-1)^{\rk(\M)} \hat{Y}_{\M}(t).
\end{align}
This modification ensures that $Y_{\M}(t)$ has non-negative coefficients, as shown in Proposition \ref{inv_z_non_negative}, rather than $\hat{Y}_\M(t)$.
For matroids with loops, 
we define 
$Y_{\M}(t):=Y_{\M \setminus \text{\{loops\}}}(t).$
Since $Z$-polynomials and inverse $Z$-polynomials of matroids are defined in a similar manner, it is natural to ask whether they share common properties.
In addition to the non-negativity of the coefficients, we also establish that the inverse $Z$-polynomials of matroids are palindromic, as shown in Lemma \ref{palindromic_inv_z}, and multiplicative under direct sums, as shown in Proposition \ref{direct_sum},
which mirror the well-known properties of the $Z$-polynomials of matroids.

In this paper, we first focus on computing explicit expressions for the inverse $Z$-polynomials of uniform matroids and sparse paving matroids.
Let $\U_{k,n}$ denote the uniform matroid of rank $k$ on $n$ elements, where the ground set is $E = [n]$ 
and the set of bases is $\mathscr{B} = \binom{[n]}{k}$. 
The first main result is the following theorem.

\begin{theorem}\label{main-theorem-uniform}
	For any uniform matroid $\U_{k,n}$ with $n \geq k\geq 1$, we have
	$$Y_{\U_{k,n}}(t)=
 \sum_{i=0}^{\lfloor k/2 \rfloor}\binom{n}{i}\binom{n-i-1}{n-k}t^i+\sum_{i=0}^{\lfloor (k-1)/2 \rfloor}\binom{n}{i}\binom{n-i-1}{n-k}t^{k-i}.$$
\end{theorem}

To compute the inverse $Z$-polynomials of sparse paving matroids, we first prove that the inverse $Z$-polynomial is a valuative invariant.

\begin{theorem}\label{z-inverse-valuative}
The assignment $\M\rightarrow Y_\M(t)$ is a valuative invariant.
\end{theorem}

The explicit formula for the inverse $Z$-polynomials of sparse paving matroids is given by the following expression.

\begin{theorem}\label{sparsing_paving_matroid-inverse-Z}
Let $\M$ be a sparse paving matroid of rank $k$ and cardinality $n$. Suppose $\M$ has $\lambda$ circuit-hyperplanes. Then
$$ 
 Y_{\M}(t) = Y_{\U_{k,n}}(t) - \lambda \left((1+t)^k - \frac{1}{2}\left((-1)^k
 +1\right)C_{\frac{k}{2}}t^{\frac{k}{2}}\right),
$$
where $C_n = \frac{1}{n+1} \binom{2n}{n}$ is the $n$-th Catalan number.
\end{theorem}

The second part of this paper explores the unimodality and log-concavity of inverse $Z$-polynomials. 
Recall that a polynomial $f(t) = a_0 + a_1 t + \cdots + a_n t^n \in \mathbb{R}[t]$ is unimodal if there exists an index $0 \leq i \leq n$ such that
$$a_0 \leq a_1 \leq \cdots \leq a_i \geq a_{i+1} \geq \cdots \geq a_n.$$
A polynomial is log-concave if
$$a_i^2 \geq a_{i-1} a_{i+1} $$ for all $0 < i < n,$
and it has no internal zeros if there do not exist indices $0 \leq i < j < k \leq n$ such that $a_i, a_k \neq 0$ and $a_j = 0$. 
Additionally, a polynomial is said to be real-rooted if all of its roots are real, or if $f(t) = 0$.
For a polynomial with nonnegative coefficients, if it has only real roots, then its coefficients are log-concave and have no internal zeros; in particular, they form a unimodal sequence.
Proudfoot, Xu, and Young \cite{PXY_2018ELC} conjectured that for any matroid $\M$, the roots of its $Z$-polynomial $Z_{\M}(t)$ lie exclusively on the negative real axis. However, computational experiments show that the inverse $Z$-polynomials may not always have real roots.
We make the following conjectures for the inverse $Z$-polynomials of matroids.

\begin{conjecture}\label{conj-uni}
\emph{For any matroid $\M$, the coefficients of the inverse $Z$-polynomial $Y_{\M}(t)$ are unimodal.}
\end{conjecture}

\begin{conjecture}\label{conj-log}
\emph{For any matroid $\M$, the coefficients of the inverse $Z$-polynomial $Y_{\M}(t)$ form a log-concave sequence with no internal zeros.
}
\end{conjecture}

Based on Theorem \ref{sparsing_paving_matroid-inverse-Z}, we prove that Conjecture \ref{conj-uni} and \ref{conj-log} hold for sparse paving matroids, as stated in the following theorem.

\begin{theorem}\label{hjkl-1}
Let $\M$ be a sparse paving matroid of rank $k$ and cardinality $n$. Then
\begin{itemize}
\item[$(1)$] the coefficients of the inverse $Z$-polynomial $Y_{\M}(t)$ are unimodal;
\item[$(2)$] the coefficients of the inverse $Z$-polynomial $Y_{\M}(t)$ form a log-concave sequence with no internal zeros.
\end{itemize}
\end{theorem}

This paper is organized as follows. In Section \ref{Matroid_Inverse_Z-polynomials}, we establish several properties of the inverse $Z$-polynomials of matroids, including their non-negativity and multiplicativity. In Section \ref{Matroid_polytopes_and_subdivisions}, we first review the basic definitions of matroid polytopes, subdivisions, and valuations. Additionally, we demonstrate that the inverse $Z$-polynomial of a matroid is a valuative invariant. In Section \ref{uni_inv_z}, we provide an explicit expression for the inverse $Z$-polynomials of uniform matroids. In Section \ref{el_sp_ma}, we derive a general formula for the inverse $Z$-polynomials of sparse paving matroids. In Section \ref{4}, we prove the unimodality and log-concavity of the inverse $Z$-polynomials of sparse paving matroids.

\section{Inverse $Z$-polynomials of matroids}\label{Matroid_Inverse_Z-polynomials}

In this section, we study several key properties of the inverse $Z$-polynomials of matroids, analogous to those of matroid $Z$-polynomials. 

We first establish the non-negativity of the coefficients of $Z$-polynomials. Braden, Huh, Matherne, Proudfoot, and Wang \cite{BHM_ARX} showed that the coefficients of matroid $Z$-polynomials are non-negative. A similar result holds for inverse $Z$-polynomials:

\begin{proposition}\label{inv_z_non_negative}
For any matroid $\M$, the coefficients of the inverse $Z$-polynomial $Y_{\M}(t)$ are non-negative.
\end{proposition}

\proof
Let $\M$ be a loopless matroid. Recall from Proudfoot \cite{proudfoot2018algebraic} that its characteristic polynomial is defined as
\begin{align}\label{char}
\chi_{\M}(t) = \sum_{F \in \mathcal{L}(\M)} \mu_{\M|_F} t^{\rk(\M/F)}.
\end{align}
The polynomial $\chi_{\M}(t)$ is monic of degree $\rk(\M)$, and its constant term is the Möbius invariant $\mu_{\M}$. Considering the alternating signs in $\chi_{\M}(t)$ (see \cite[Theorem 4]{Rota}), we deduce that $\mu_{\M}$ is $ (-1)^{\rk(\M)}$ multiplied by a positive integer.

Since the inverse Kazhdan–Lusztig polynomial $Q_{\M}(t)$ has non-negative coefficients, it follows directly from \eqref{z-inverse-algebra} that each summation term
$$ (-1)^{\rk(\M|_F)} Q_{\M|_F}(t)\, t^{\rk(\M/F)} \mu_{\M/F} $$
in the definition of $\hat{Y}_{\M}(t)$ is $(-1)^{\rk(\M)}$ times a polynomial with non-negative coefficients. Thus, the polynomial 
$Y_{\M}(t)=(-1)^{\rk(\M)}\hat{Y}_{\M}(t)$ 
has non-negative coefficients.
\qed

With non-negativity established, we turn to palindromicity. It is known that matroid $Z$-polynomials are palindromic \cite{PXY_2018ELC}. Braden, Huh, Matherne, Proudfoot, and Wang \cite[p.5]{BHM_ARX} observed that inverse $Z$-polynomials also share this property; thus, we state the result without proof.

\begin{lemma}\label{palindromic_inv_z}
 For any matroid $\M$, the polynomial $Y_\M(t)$ is palindromic of degree $\rk(\M)$, i.e.,
 $$
 t^{\rk(\M)} Y_\M(t^{-1}) = Y_\M(t).
 $$
\end{lemma}

We now demonstrate that the inverse $Z$-polynomial of a matroid is multiplicative with respect to direct sums. The direct sum of two matroids $\M_1 = (E_1, \mathscr{B}_1)$ and $\M_2 = (E_2, \mathscr{B}_2)$, where $E_1 \cap E_2 = \emptyset$, is defined on the ground set $E_1 \cup E_2$ with the set of bases
$$ \mathscr{B}_1 \cup \mathscr{B}_2 = \{ B_1 \cup B_2 : B_1 \in \mathscr{B}_1, B_2 \in \mathscr{B}_2 \}. $$ 
It is a well-known fact that the (inverse) Kazhdan-Lusztig polynomials and the $Z$-polynomials of matroids are all multiplicative under direct sums; see \cite
{elias2016kazhdan}, \cite
{gao2021inverse}, and \cite{PXY_2018ELC}, respectively. 
The corresponding result for inverse $Z$-polynomials is as follows.

\begin{proposition}\label{direct_sum}
For any matroids $\M_1$ and $\M_2$ with respective ground sets $E_1$ and $E_2$, we have
$$Y_{\M_1\oplus \M_2}(t)=Y_{\M_1}(t) \cdot Y_{ \M_2}(t).$$
\end{proposition}

\proof
The statement can be proven in the same manner as Elias, Proudfoot and Wakefield did for the matroid Kazhdan-Lusztig polynomials. From \eqref{defi-z-mod}, 
it suffices to show that
$$\Hat{Y}_{\M_1\oplus \M_2}(t)=\Hat{Y}_{\M_1}(t)\Hat{Y}_{ \M_2}(t).$$
We proceed by induction on the sum $\mathrm{rk} (\M_1)+\mathrm{rk}(\M_2)$. The statement is clear when $\mathrm{rk} (\M_1)=0$ or $\mathrm{rk}(\M_2)=0$. Now, assume that the statement holds for $\M_1'$ and $\M_2'$ whenever one of $\mathrm{rk} (\M_1')\leq \mathrm{rk} (\M_1)$ and $\mathrm{rk} (\M_2')\leq \mathrm{rk}(\M_2)$, with strict inequalities in at least one of the ranks.

By the properties of flats, we have
$$\mathcal{L}(\M_1\oplus \M_2)=\mathcal{L}(\M_1)\times \mathcal{L}(\M_2)~~~~~~
\text{and}
~~~~~~\mathrm{rk} \, (\M_1\oplus \M_2)=\mathrm{rk}\, (\M_1)+\mathrm{rk}\, (\M_2).$$
Moreover, for any $(F_1,F_2)\in \mathcal{L}(\M_1\oplus \M_2)$,
\begin{align*}
 (\M_1\oplus \M_2)|_{(F_1,F_2)} & = (M_1|_{F_1}) \oplus (\M_2|_{F_2}),\\[6pt]
 (\M_1\oplus \M_2)/{(F_1,F_2)} & = (\M_1/{F_1}) \oplus (\M_2/{F_2}).
\end{align*}
Hence, using \eqref{z-inverse-algebra}, we obtain
\begin{align*}
	\Hat{Y}_{\M_1\oplus \M_2}(t)
 &=\sum_{(F_1,F_2)\in \mathcal{L}(\M_1\oplus \M_2)}
 \Hat{Q}_{(\M_1\oplus \M_2)|_{(F_1,F_2)}} (t) \cdot t^{\mathrm{rk}\left( (\M_1\oplus \M_2)/ {(F_1,F_2) } \right)} \mu_{(\M_1\oplus \M_2)/{(F_1,F_2)}}\\[6pt]
 &=\sum_{(F_1,F_2)\in \mathcal{L}(\M_1\oplus \M_2)}\Hat{Q}_{(\M_1|_{F_1}) \oplus (\M_2|_{F_2})}(t) \cdot t^{\mathrm{rk} \left( \M_1/{F_1} \right)+\mathrm{rk}\left( \M_2/{F_2} \right)} \mu_{(\M_1/{F_1}) \oplus (\M_2/{F_2})}.
\end{align*}
By the inductive hypothesis, we get
\begin{align}
\Hat{Y}_{\M_1\oplus \M_2}(t)
&=
\Hat{Q}_{(\M_1|_{E_1}) \oplus (\M_2|_{E_2})}(t) \cdot t^{\mathrm{rk} \left(\M_1/{E_1}\right)+\mathrm{rk} \left( \M_2/{E_2} \right) } \mu_{(\M_1/{E_1}) \oplus (\M_2/{E_2})}
\nonumber\\[5pt]
&\quad
+\sum_{(F_1,F_2)\neq (E_1, E_2)}\Hat{Q}_{(\M_1|_{F_1} )\oplus (\M_2|_{F_2})}(t) \cdot t^{\mathrm{rk}\left( \M_1/{F_1} \right)+\mathrm{rk}\left( \M_2/{F_2}\right)}\mu_{(\M_1/{F_1} )\oplus (\M_2/{F_2})}\nonumber\\[6pt]
&=\Hat{Q}_{\M_1 \oplus \M_2}(t)\nonumber\\[5pt]
&\quad
+\sum_{(F_1,F_2)\neq (E_1, E_2)}\left(\hat{Q}_{\M_1|_{F_1}}(t) \cdot t^{\mathrm{rk} \left( \M_1/{F_1} \right)} \mu_{\M_1/{F_1}}\right)\cdot \left(\hat{Q}_{\M_2|_{F_2}}(t) \cdot t^{\mathrm{rk} \left( \M_2/{F_2} \right)} \mu_{\M_2/{F_2}}\right),\label{eq-before}
\end{align}
where the last identity follows from the fact that the M\"{o}bius invariants of matroids are multiplicative under direct sums, as shown in \cite[Theorem 7.1.7.]{Zaslavsky1987}.
On the other hand, by \eqref{z-inverse-algebra}, we have
\begin{align*}
\Hat{Y}_{\M_1}(t)&= \sum_{F_1\in \mathcal{L}(\M_1)}
 \Hat{Q}_{\M_1|_{F_1}}(t) \cdot t^{\mathrm{rk} (\M_1/{F_1})} \mu_{\M_1/{F_1}},\\
\Hat{Y}_{\M_2}(t)&= \sum_{F_2\in \mathcal{L}(\M_2)}
\Hat{Q}_{\M_2|_{F_2}}(t) \cdot t^{\mathrm{rk} ( \M_2/{F_2})} \mu_{\M_2/{F_2}},
\end{align*}
which leads to
\begin{align}
&\Hat{Y}_{\M_1}(t) \cdot \Hat{Y}_{\M_2}(t)
=\Hat{Q}_{\M_1}(t)\cdot \Hat{Q}_{\M_2}(t)\nonumber\\[5pt]
&\quad
+\sum_{(F_1,F_2)\neq (E_1, E_2)}
\left(\Hat{Q}_{\M_1|_{F_1}}(t) \cdot t^{\mathrm{rk} (\M_1/{F_1})} \mu_{\M_1/{F_1}}\right)\cdot 
\left(\Hat{Q}_{\M_2|_{F_2}}(t) \cdot t^{\mathrm{rk} (\M_2/{F_2})} \mu_{\M_2/{F_2}}\right).\label{eq-after}
\end{align}
Subtracting \eqref{eq-after} from \eqref{eq-before} yields
\begin{align*}
\Hat{Y}_{\M_1\oplus \M_2}(t)-\Hat{Y}_{\M_1}(t)\Hat{Y}_{\M_2}(t)=	0.
\end{align*}
This completes the proof.
\qed

\section{Matroid polytopes, subdivisions, and valuations }\label{Matroid_polytopes_and_subdivisions}

The aim of this section is to prove Theorem \ref{z-inverse-valuative}. To this end, we first review the basic definitions of matroid polytopes, subdivisions, and valuations, primarily following the notions and notation summarized by Ferroni and Schr\"{o}ter \cite{ferroni2022valuative}.

Let $\M=(E,\mathscr{B})$ be a matroid. For each $i \in E$, let $e_i$ denote the $i$-th canonical vector in $\mathbb{R}^E$. The base polytope of $\M$, denoted by $\PP(\M)$, is defined as the convex hull of the set of vectors $e_B$, where each $B \in \mathscr{B}$ is a base of $\M$, and $e_B = \sum_{i \in B} e_i$. That is,
$$ \PP(\M) := \text{conv}\{ e_B : B \in \mathscr{B} \} \subseteq \mathbb{R}^E. $$
A subdivision of $\PP = \PP(\M) \subseteq \mathbb{R}^E$ is a finite collection 
$\mathcal{S} = \{\PP_1, \PP_2, \ldots, \PP_s\}$ of matroid polytopes $\PP_i = \PP(\M_i) \subseteq \mathbb{R}^E$, such that the following conditions hold:
\begin{itemize}
 \item $\PP = \bigcup_{i=1}^s \PP_i$.
 \item Any face of a polytope $\PP_i \in \mathcal{S}$ is contained in some polytope in $\mathcal{S}$.
 \item For each pair of distinct indices $i \neq j$, the intersection $\PP_i \cap \PP_j$ is a common (possibly empty) face of both $\PP_i$ and $\PP_j$.
\end{itemize}
Sometimes, when dealing with a subdivision $\mathcal{S}$ of a matroid polytope $\PP$, it is convenient to focus only on the interior faces of the subdivision. More precisely, if we denote the boundary of $\PP$ by $\partial \PP$, we define
$$ \mathcal{S}^{\intt} := \{\PP_i \in \mathcal{S} : \PP_i \not \subseteq \partial \PP \}. $$
Throughout this paper, we consider matroids with ground sets of the form $E \subseteq \mathbb{Z}_{\geq 1}$, where $|E| < \infty$. For any such $E$, let $\mathcal{M}_E$ denote the set of all base polytopes of matroids on $E$, and let $\mathcal{M}$ be the union of all the $\mathcal{M}_E$'s.
Let $A$ be an abelian group (typically a ring). 
Consider a map $f: \mathcal{M} \to A$. If $f\left(\mathscr{P}(\M)\right) = f\left(\mathscr{P}(\mathrm{N})\right)$ whenever $\M \cong \mathrm{N}$, we say that $f$ is an \emph{invariant}. 
If $f$ satisfies the following property for every matroid polytope $\mathscr{P}$ and every subdivision $\mathcal{S} = \{\mathscr{P}_1, \mathscr{P}_2, \ldots, \mathscr{P}_s\}$ of $\mathscr{P}$, i.e.,
$$ f(\mathscr{P}) = \sum_{\mathscr{P}_i \in \mathcal{S}^{\mathrm{int}}} (-1)^{\dim \mathscr{P} - \dim \mathscr{P}_i} f(\mathscr{P}_i), $$
we say that $f$ is a \emph{valuation}.
To simplify our notation, we will use $f(\M)$ to denote the map $f: \mathcal{M} \to A$ instead of writing $f(\PP(\M))$.

In order to prove Theorem \ref{z-inverse-valuative}, we first need to establish the following basic lemma on matroid valuative invariants, which has appeared in several previous works. For the sake of completeness, we provide a self-contained proof.

\begin{lemma}\label{lemma-invaluative}
If $f: \mathcal{M} \rightarrow \mathbb{R}[t]$ is a valuative invariant, then
the maps $g_1,g_2: \mathcal{M} \rightarrow \mathbb{R}[t]$ defined by
\begin{align*}
g_1(\M)=t^{\rk (\M)}f(\M)~~~~~~~~~~\text{and}~~~~~~~~~~
g_2(\M)=(-1)^{\rk (\M)}f(\M),
\end{align*}
are both valuative invariants
\end{lemma}

\proof 
Since $f: \mathcal{M} \rightarrow \mathbb{R}[t]$ is a valuative invariant, we have
$$
f\left(\M\right)=\sum_{\mathscr{P}(\M_i) \in \mathcal{S}^{\mathrm{int}} }(-1)^{\dim \mathscr{P}(\M)-\dim \mathscr{P}(\M_i) }f(\M_i)
$$
for any matroid $\M$ and any subdivision $\mathcal{S}=\{\mathscr{P}(\M_1) ,\mathscr{P}(\M_2) ,\ldots,\mathscr{P}(\M_s)   \}$ of $\PP=\PP(\M)$ in $ \mathcal{M}$.
Note that all matroids corresponding to the polytopes in $\mathcal{S}$ share the same ground set and rank. Therefore, we have
\[
t^{\rk(\M)} = t^{\rk(\M_i)} \quad \text{and} \quad (-1)^{\rk(\M)} = (-1)^{\rk(\M_i)}
\]
for each $\M_i$. The desired result then follows, which completes the proof.
\qed

Let $R$ be a ring and let $f, g: \mathcal{M} \to R$ be two maps. The convolution of $f$ and $g$, denoted by $f \star g: \mathcal{M} \to R$, is defined as
$$
(f \star g)(\M) := \sum_{S \subseteq E} f(\M|_S) g(\M/S),
$$
where $\M$ is a matroid with ground set $E$.
It is worth noting that if $g: \mathcal{M} \to R$ satisfies $g(\M) = 0$ whenever $\M$ contains loops, the convolution of any other invariant $f$ with $g$ depends only on the flats of $\M$. Specifically, we have
$$(f\star g)(\M)=\sum_{F\in \mathcal{L}(\M)} f(\M|_F)g(\M/F),$$
since the contraction $\M/S$ is loopless if and only if $S$ is a flat.
One of the main results of Ardila and Sanchez \cite{Ardila_2023valuations} establishes that the convolution of two valuations
is itself a valuation, as stated in the following theorem.

\begin{theorem}[\cite{Ardila_2023valuations}, Theorem C]\label{f_g_convolution}
Let $R$ be a ring. If $f, g: \mathcal{M} \rightarrow A$ are valuations, then the convolution $f \star g$ is
also a valuation.
\end{theorem}

We are now ready to prove that the inverse $Z$-polynomial is a valuative invariant.

\proof[Proof of Theorem \ref{z-inverse-valuative}]
Let $f,g: \mathcal{M} \rightarrow \mathbb{R}[t]$ be defined by
$$
f(\M)= (-1)^{\mathrm{rk} (\M)} Q_{\M}(t) \quad \text{and} \quad g(\M) = t^{\mathrm{rk} (\M)} \mu_{\M}.
$$
It is known from Ardila and Sanchez \cite[Theorem 8.8]{Ardila_2023valuations} that the map $\M \rightarrow Q_{\M}(t)$ is a valuative invariant. By applying Lemma \ref{lemma-invaluative}, it follows that $f$ is also a valuative invariant.
Furthermore, recall that the Tutte polynomial $T_{\M}(x,y)$ of a matroid $\M=(E,\mathscr{B})$ is defined as
$$
T_{\M}(x,y) := \sum_{A \subseteq E} (x-1)^{\mathrm{rk} (E) - \mathrm{rk} (A)} (y-1)^{|A| - \mathrm{rk} (A)} \in \mathbb{Z}[x,y].
$$
There are several established proofs showing that the map $ \mathcal{M} \rightarrow \mathbb{Z}[x,y]$ defined by $f: \M \rightarrow T_{\M}(x,y)$ is a valuative invariant; see \cite[Theorem 7.17]{ferroni2022valuative}.
From the identity 
$$
\chi_{\M}(t) = (-1)^{\mathrm{rk}(\M)} T_{\M}(1-t,0),
$$
we deduce that both the characteristic polynomial $\chi_{\M}(t)$ and the M\"{o}bius invariant $\mu_{\M} = \chi_{\M}(0)$ are valuative invariants. Again, by Lemma \ref{lemma-invaluative}, we conclude that $g$ is a valuative invariant as well.
It is worth noting that $g = t^{\mathrm{rk} (\M)} \mu_{\M} = 0$ whenever $\M$ has loops. Combining this with Theorem \ref{f_g_convolution}, we obtain the following expression for the convolution
\begin{align*}
(f \star g)(\M) &= \sum_{F \in \mathcal{L}(\M)} f(\M|_F) g(\M/F) \\
&= \sum_{F \in \mathcal{L}(\M)} (-1)^{\mathrm{rk}(F)} Q_{\M|_F}(t) \cdot t^{\mathrm{rk}(\M/F)} \mu_{\M/F} \\
&= (-1)^{\mathrm{rk}(\M)} Y_{\M}(t),
\end{align*}
which is a valuative invariant. This completes the proof.
\qed

\section{Uniform matroids}\label{uni_inv_z}

The aim of this section is to compute the inverse $Z$-polynomials of uniform
matroids. 
To this end, 
we first provide the expression for $\mu_{\M}$ in the case of uniform matroids.
\begin{lemma}[{\cite[p.121]{Zaslavsky1987}}]
 For any uniform matroid $\U_{k,n}$ with $n \geq k \geq 1$, its M\"{o}bius invariant satisfies
 \begin{align}\label{mobius_uniform}
 \mu_{\U_{k,n}}= (-1)^k \binom{n-1}{k-1}.
 \end{align}
\end{lemma}

We are now ready to prove Theorem \ref{main-theorem-uniform}, which presents an explicit formula for the inverse $Z$-polynomials of uniform matroids.

\proof[Proof of Theorem \ref{main-theorem-uniform}]
From Lemma \ref{palindromic_inv_z}, we know that for any matroid $\M$, the polynomial $Y_\M(t)$ is palindromic of degree $\rk (\M)$. Therefore, it suffices to show that for any $0 \leq i \leq \lfloor k / 2 \rfloor$, the coefficient of $t^i$ in $Y_{\U_{k,n}}(t)$ is given by
$
\binom{n}{i}\binom{n-i-1}{n-k}.
$

Applying equation \eqref{z-inverse-algebra} to uniform matroids, we obtain
$$
\hat{Y}_{\U_{k,n}}(t) = \sum_{F \in \mathcal{L}(\U_{k,n})} (-1)^{\mathrm{rk} (F)} Q_{\U_{k,n}|_F}(t) \cdot t^{\mathrm{rk}(\U_{k,n}) - \mathrm{rk}(F)} \mu_{\U_{k,n}/F}.
$$
Given $0 \leq i <k$, let $F$ denote any flat of $\U_{k,n}$ with rank $i$. It is clear that the restriction $\U_{k,n}|_F \cong B_i$ and the contraction $\U_{k,n}/F \cong \U_{k-i,n-i}$. 
Here, $B_n$ denotes the Boolean matroid, a special uniform matroid of rank $n$ and cardinality $n$ for $n \geq 1$, and $B_0$ denotes the unique matroid of rank 0.
Moreover, there exist $\binom{n}{i}$ such flats. For $i = k$, there exists only one flat of rank $k$, which is the ground set $E$. Hence, we have
\begin{align}\label{equa-inverse-Z-uniform-1}
\hat{Y}_{\U_{k,n}}(t) = \sum_{i=0}^{k-1} \binom{n}{i} (-1)^i Q_{B_i}(t) \cdot t^{k-i} \mu_{\U_{k-i,n-i}} + (-1)^k Q_{\U_{k,n}}(t).
\end{align}
From \cite[Corollary 3.2]{gao2021inverse}, we immediately obtain the result
\begin{align}\label{equa-inverse-KL-uniform-1}
Q_{B_n}(t) = 1
\end{align}
for any $n \geq 1$.
Additionally, it is known that
\begin{align}\label{equa-inverse-KL-uniform-2}
Q_{\U_{k,n}}(t) = \binom{n}{k} \sum_{i=0}^{\lfloor (k-1) / 2 \rfloor} \frac{(n-k)(k-2i)}{(n-k+i)(n-i)} \binom{k}{i} t^i
\end{align}
for $n > k \geq 1$; see \cite[Theorem 3.3]{gao2021inverse}. 
Substituting equations \eqref{mobius_uniform}, \eqref{equa-inverse-KL-uniform-1}, and \eqref{equa-inverse-KL-uniform-2} into \eqref{equa-inverse-Z-uniform-1}, we obtain
\begin{align*}
\hat{Y}_{\U_{k,n}}(t) &= \sum_{i=0}^{k-1} \binom{n}{i} (-1)^i \cdot (-1)^{k-i} \binom{n-i-1}{k-i-1} t^{k-i} \\
&\quad + (-1)^k \binom{n}{k} \sum_{i=0}^{\lfloor (k-1) / 2 \rfloor} \frac{(n-k)(k-2i)}{(n-k+i)(n-i)} \binom{k}{i} t^i \\
&= (-1)^k \left( \sum_{i=0}^{k-1} \binom{n}{i} \binom{n-i-1}{n-k} t^{k-i} + \binom{n}{k} \sum_{i=0}^{\lfloor (k-1) / 2 \rfloor} \frac{(n-k)(k-2i)}{(n-k+i)(n-i)} \binom{k}{i} t^i \right).
\end{align*}
Thus, we obtain
\begin{align}
Y_{\U_{k,n}}(t) &= \sum_{i=0}^{k-1} \binom{n}{i} \binom{n-i-1}{n-k} t^{k-i} + \binom{n}{k} \sum_{i=0}^{\lfloor (k-1) / 2 \rfloor} \frac{(n-k)(k-2i)}{(n-k+i)(n-i)} \binom{k}{i} t^i.
\end{align}
Since $Y_{\U_{k,n}}(t)$ is palindromic of degree $k$, we have
\begin{align*}
Y_{\U_{k,n}}(t) &= t^k Y_{\U_{k,n}}(t^{-1}) \\
&= \sum_{i=0}^{k-1} \binom{n}{i} \binom{n-i-1}{n-k} t^i + \binom{n}{k} \sum_{i=0}^{\lfloor (k-1) / 2 \rfloor} \frac{(n-k)(k-2i)}{(n-k+i)(n-i)} \binom{k}{i} t^{k-i}.
\end{align*}
Note that for any positive integer $k$, we have $\lfloor k / 2 \rfloor + \lfloor (k-1) / 2 \rfloor = k-1 < k$, which implies that
$$
\lfloor k / 2 \rfloor < k - \lfloor (k-1) / 2 \rfloor.
$$
Therefore, for any $0 \leq i \leq \lfloor k / 2 \rfloor$, the coefficient of $t^i$ in $Y_{\U_{k,n}}(t)$ is
$\binom{n}{i} \binom{n-i-1}{n-k}.$
Thus, we complete the proof.
\qed

Note that the proof of Theorem \ref{main-theorem-uniform} only relies on the evaluation of the M\"{o}bius invariants and inverse Kazhdan-Lusztig polynomials of uniform matroids. Once Theorem \ref{main-theorem-uniform} is established, it becomes straightforward to compute the $Z$-polynomials of uniform matroids.
Several formulas for $Z_{\U_{m,d}}(t)$ have been obtained in the literature. For instance, see Proudfoot, Xu, and Young \cite{PXY_2018ELC}, and Gao, Lu, Xie, Yang, and Zhang \cite{gao2021uniform}. The following result presents a new formula for $Z_{\U_{m,d}}(t)$.

\begin{corollary}
	For any uniform matroid $\U_{k,n}$ with $n \geq k\geq 1$, we have
	\begin{align*}
 Z_{\U_{k,n}}(t)
 = & \sum_{j=0}^{\lfloor k/2 \rfloor} \sum_{i=0}^{k-2j} \frac{(-1)^{k-i+1}}{n-i-j} \binom{n}{1,i,j,n-k,k-i-j-1} (1+t)^i t^j \\
 &\quad + \sum_{j=0}^{\lfloor (k-1)/2 \rfloor} \sum_{i=0}^{k-2j-1} \frac{(-1)^{k-i+1}}{n-i-j} \binom{n}{1,i,j,n-k,k-i-j-1} (1+t)^i t^{k-i-j},
\end{align*}
where $\binom{n}{a_1, a_2, \ldots,a_k}=\frac{n!}{a_1!a_2!\cdots a_k!}$ denotes the multinomial coefficient.
\end{corollary}

\proof
From the general theory of Kazhdan-Lusztig-Stanley polynomials on $\mathcal{L}(\U_{k,n})$, for $n \geq k \geq 1$, we have
\begin{align*}
Z_{\U_{k,n}}(t)
=& -\sum_{E \neq F \in \mathcal{L}(\U_{k,n})}
Z_{\U_{k,n}|_F}(t) \cdot (-1)^{\rk(\U_{k,n}/F)} Y_{\U_{k,n}/F}(t) \\
=& -\sum_{i=0}^{k-1} \binom{n}{i} Z_{B_i}(t) \cdot (-1)^{\rk(\U_{k-i,n-i})} Y_{\U_{k-i,n-i}}(t).
\end{align*}
For $0 \leq i \leq k-1$, we have $P_{B_i}(t) = 1$ (see \cite[Proposition 2.7]{elias2016kazhdan}), and using \eqref{z_defini} applied to $B_i$, we get
$Z_{B_i}(t) = (1+t)^i.$
From Theorem \ref{main-theorem-uniform}, we know
\begin{align*}
 Y_{\U_{k-i,n-i}}(t) = \sum_{j=0}^{\lfloor (k-i)/2 \rfloor} \binom{n-i}{j} \binom{n-i-j-1}{n-k} t^j 
 + \sum_{j=0}^{\lfloor (k-i-1)/2 \rfloor} \binom{n-i}{j} \binom{n-i-j-1}{n-k} t^{k-i-j}.
\end{align*}
Substituting these into the expression for $Z_{\U_{k,n}}(t)$, we get
\begin{align*}
 Z_{\U_{k,n}}(t)
 &= -\sum_{i=0}^{k-1} \binom{n}{i} (1+t)^i
 \cdot (-1)^{k-i} \Big[ \sum_{j=0}^{\lfloor (k-i)/2 \rfloor} \binom{n-i}{j} \binom{n-i-j-1}{n-k} t^j \\
 &\quad + \sum_{j=0}^{\lfloor (k-i-1)/2 \rfloor} \binom{n-i}{j} \binom{n-i-j-1}{n-k} t^{k-i-j} \Big] \\
 &= \sum_{j=0}^{\lfloor k/2 \rfloor} \sum_{i=0}^{k-2j} \frac{(-1)^{k-i+1}}{n-i-j} \binom{n}{1,i,j,n-k,k-i-j-1} (1+t)^i t^j \\
 &\quad + \sum_{j=0}^{\lfloor (k-1)/2 \rfloor} \sum_{i=0}^{k-2j-1} \frac{(-1)^{k-i+1}}{n-i-j} \binom{n}{1,i,j,n-k,k-i-j-1} (1+t)^i t^{k-i-j},
\end{align*}
where the last equality follows from interchanging the order of summation. This completes the proof.
\qed

\section{Sparse paving matroids}\label{el_sp_ma}

This section is dedicated to proving Theorem \ref{sparsing_paving_matroid-inverse-Z}. In fact, based on Theorem \ref{z-inverse-valuative}, we derive a general formula for the inverse $Z$-polynomials of large classes of matroids, including elementary split matroids and (sparse) paving matroids.

Let us first review the definitions of elementary split matroids and (sparse) paving matroids.
Let $\M=(E,\mathscr{B})$ be a matroid. An independent set of $\M$ is a set $I \subseteq E$ contained in some $B \subseteq \mathscr{B}$. If a set is not independent, we will say that it is dependent. A circuit of $\M$ is a minimal dependent set.
A matroid $\M$ is called \emph{paving} if all its circuits have size at least $\rk(\M)$. If both $\M$ and its dual $\M^*$ are paving, then $\M$ is referred to as a \emph{sparse paving} matroid. A well-known conjecture by Mayhew, Newman, Welsh, and Whittle \cite{mayhew2011asymptotic} affirms that asymptotically, almost all matroids are sparse paving.
The class of elementary split matroids is somewhat more complicated to define. First introduced in \cite{BKSYY_2023}, they form a subclass of the split matroids introduced by Joswig and Schröter in \cite{JS_2017}. Ferroni and Schröter \cite{ferroni2022valuative} conjectured that the class of elementary split matroids is strictly larger than the class of sparse paving matroids. It is expected that elementary split matroids will dominate even when considering only non sparse paving matroids.

The aforementioned classes of matroids can be equivalently defined through the relaxation of stressed subsets.
Let $\M$ be a matroid on the ground set $E$ with rank $k$. A subset $A \subseteq E$ is called stressed if both the restriction $\M|_A$ and the contraction $\M/A$ are isomorphic to uniform matroids. Ferroni and Schröter \cite[Section 3]{ferroni2022valuative} describe how to construct a new matroid on $E$ by relaxing a stressed subset $A$, which generalizes the well-known circuit-hyperplane relaxation (A flat of rank $\rk(\M)-1$ is said to be a hyperplane).
Precisely, the cusp of $A$ is defined as the collection of $k$-sets
$$ \cusp_{\M}(A) := \left\{ S \in \binom{E}{k} : |S \cap A| \geq \rk(A) + 1 \right\}, $$
and a new matroid is constructed on the ground set $E$ whose family of bases is given by $\mathscr{B} \cup \cusp_{\M}(A)$. This matroid is referred to as the relaxation of $\M$ by $A$ and is denoted by $\Rel(\M, A)$.
If $A$ and $A'$ are two distinct stressed subsets of $\M$, then $A'$ remains stressed in $\Rel(\M, A)$, and the cusp of $A'$ in $\Rel(\M, A)$ is given by $\cusp_{\M}(A')$. 
In particular, for each matroid, one may associate a canonical relaxed matroid, obtained by relaxing (in any order) all stressed subsets with non-empty cusp. 
The following proposition provides equivalent definitions of the previously discussed classes of matroids in terms of the relaxation of stressed subsets. For further details, we refer to \cite[Corollary 3.17]{ferroni2023stressed}, \cite[Theorem 4.8]{ferroni2022valuative}, and
\cite[Lemma 5.1]{ferroni2022matroid}.

\begin{proposition} \label{s-paving-equivalent-proposition}
 Let $\M$ be a matroid of rank $k$ on the ground set $E=[n]$.
 \begin{itemize}
 \item[$(1)$] $\M$ is sparse paving if and only if the relaxation of all of its circuit-hyperplanes yields the matroid $\U_{k,n}$.
 \item[$(2)$] $\M$ is paving if and only if the relaxation of all of its stressed hyperplanes of size at least $k$ yields the matroid $\U_{k,n}$.
 \item[$(3)$] $\M$ is elementary split if and only if the relaxation of all of its stressed subsets with non-empty cusp yields the matroid $\U_{k,n}$.
 \end{itemize}
\end{proposition}

In a direct sum of two uniform matroids, the ground sets of each of the summands
are stressed and therefore can be relaxed. 
For $0 \leq r \leq h$ and $0\leq k - r \leq n - h$, Ferroni and Schr\"{o}ter \cite{ferroni2022valuative} defined the cuspidal matroid 
$\LL_{r,k,h,n} := \Rel(\U_{k-r,n-h}\oplus\U_{r,h}, \U_{r,h}),$
 i.e., they relaxed the subset corresponding to the ground set of the direct summand $\U_{r,h}$.
We have the following general expression for $Y_{\M}(t)$ when $\M$ is an elementary split matroid.

\begin{proposition}\label{inver_Z_elemnt-split}
Let $\M$ be an elementary split matroid of rank $k$ and cardinality $n$. Then
 $$ 
 Y_{\M}(t) = Y_{\U_{k,n}}(t) - \sum_{r,h}\lambda_{r,h} \left(Y_{\Lambda_{r,k,h,n}} (t)- Y_{\U_{k-r,n-h}}(t)\cdot Y_{\U_{r,h}}(t)\right) 
 $$
 where $\lambda_{r,h}$ denotes the number of stressed subsets with non-empty cusp of
size $h$ and rank $r$ in $\M$.
\end{proposition}

\proof
For any valuative invariant $f$ and elementary split matroid $\M$ of rank $k$ and cardinality $n$, it is known from Ferroni and Schr\"{o}ter \cite[Theorem 5.3]{ferroni2022valuative} that 
 $$ 
 f(\M) = f(\U_{k,n}) - \sum_{r,h}\lambda_{r,h} \left(f(\mathsf{\Lambda}_{r,k,h,n}) - f(\U_{k-r,n-h}\oplus\U_{r,h})\right).
 $$
Since the inverse $Z$-polynomial is a valuative invariant by Theorem \ref{z-inverse-valuative} and is multiplicative under direct sums of matroids by Lemma \ref{direct_sum}, applying these properties to the above identity yields the desired result.
\qed

When $\M$ is a paving matroid, the expression for $Y_{\M}(t)$ in Proposition \ref{inver_Z_elemnt-split} simplifies as follows.

\begin{proposition}\label{theorem-Z_inverser_maving}
Let $\M$ be an paving matroid of rank $k$ and cardinality $n$. Then
$$ 
 Y_{\M}(t) = Y_{\U_{k,n}}(t) - \sum_{h \geq k}\lambda_{h} \left(Y_{\U_{k,h+1}} (t)- (1+t)\cdot Y_{\U_{k-1,h}}(t)\right),
 $$
where $\lambda_{h}$ denotes the number of stressed hyperplanes of size $h$ in $\M$.
\end{proposition}

\proof
In a paving matroid of rank $k$, the only stressed subsets with non-empty cusp are the stressed hyperplanes of size at least $k$. Since the rank of any hyperplane is $k-1$, we have the following expression for the inverse $Z$-polynomial of $\M$,
$$ 
 Y_{\M}(t) = Y_{\U_{k,n}}(t) - \sum_{h \geq k}\lambda_{h} \left(Y_{\Lambda_{k-1,k,h,n}}(t) - Y_{\U_{1,n-h}}(t) \cdot Y_{\U_{k-1,h}}(t)\right).
$$
From Ferroni, Nasr, and Vecchi \cite[Proposition 3.11]{ferroni2023stressed}, we know that the simplification of the cuspidal matroid $\Lambda_{k-1,k,h,n}$ is isomorphic to the uniform matroid $\U_{k,h+1}$. Since the inverse $Z$-polynomials remain invariant under taking simplifications of loopless matroids, it follows that
$$ 
 Y_{\Lambda_{k-1,k,h,n}}(t) = Y_{\U_{k,h+1}}(t).
$$
Moreover, by Theorem \ref{main-theorem-uniform}, we have
$Y_{\U_{1,n-h}}(t) = 1 + t.$
Combining these results, we obtain the desired expression for $Y_{\M}(t)$, completing the proof.
\qed

We now proceed to give a proof of Theorem \ref{sparsing_paving_matroid-inverse-Z}.

\proof[Proof of Theorem \ref{sparsing_paving_matroid-inverse-Z}]
Since $\M$ is a sparse paving matroid, the stressed hyperplanes and circuit-hyperplanes in $\M$ are exactly the same. Moreover, all stressed hyperplanes have the same size, which is equal to $k$. Thus, by Proposition \ref{theorem-Z_inverser_maving}, we obtain
\begin{align*} 
 Y_{\M}(t) = Y_{\U_{k,n}}(t) - \lambda \left(Y_{\U_{k,k+1}}(t) - (1+t)\cdot Y_{\U_{k-1,k}}(t)\right).
\end{align*}
It suffices to show that for $k \geq 1$,
\begin{align}\label{equation-theorem-sparse}
Y_{\U_{k,k+1}}(t) - (1+t) \cdot Y_{\U_{k-1,k}}(t) = (1+t)^k - \frac{1}{2} \left((-1)^k + 1 \right) C_{\frac{k}{2}} t^{\frac{k}{2}}.
\end{align}

From Lemma \ref{palindromic_inv_z}, we know that for any $k \geq 1$, $Y_{\U_{k,k+1}}(t)$ is a palindromic polynomial of degree $k$, with its center of symmetry at $\lfloor k/2 \rfloor$. Moreover, by Theorem \ref{main-theorem-uniform}, we have
\begin{align}\label{equation_inverse_uni-prod}
 Y_{\U_{k,k+1}}(t) = \sum_{i=0}^{\lfloor k/2 \rfloor} (k-i) \binom{k+1}{i} t^i + \sum_{i=0}^{\lfloor (k-1)/2 \rfloor} (k-i) \binom{k+1}{i} t^{k-i},
\end{align}
which tells us that the leading coefficient of $t^k$ in $Y_{\U_{k,k+1}}(t)$ is $k$.
Given the symmetry of $1+t$ and $Y_{\U_{k,k+1}}(t)$, it follows that the polynomials $(1+t) \cdot Y_{\U_{k-1,k}}(t)$ and $Y_{\U_{k,k+1}}(t)-(1+t) \cdot Y_{\U_{k-1,k}}(t)$ are both palindromic polynomials of degree $k$, with their center of symmetry at $\lfloor k/2 \rfloor$. Furthermore, the leading coefficients of these polynomials are $k-1$ and $1$, respectively.
Thus, we need to show that for $0 \leq i < \frac{k}{2}$, the coefficient of $t^i$ in $Y_{\U_{k,k+1}}(t) - (1+t) \cdot Y_{\U_{k-1,k}}(t)$ is
$ \binom{k}{i}, $
and for $k$ even, when $i = \frac{k}{2}$, the coefficient of $t^{\frac{k}{2}}$ is
$ \binom{k}{\frac{k}{2}} - C_{\frac{k}{2}}. $

To compute this, we first consider the case when $0 \leq i < \frac{k}{2}$. Note that $\binom{n+1}{m+1} - \binom{n}{m} = \binom{n}{m+1}$ for any $n,m \geq 0$, we obtain
\begin{align*}
[t^i] \left(Y_{\U_{k,k+1}}(t) - (1+t)\cdot Y_{\U_{k-1,k}}(t)\right) 
=& [t^i] Y_{\U_{k,k+1}}(t) - [t^i] Y_{\U_{k-1,k}}(t) - [t^{i-1}] Y_{\U_{k-1,k}}(t) \\
=& (k-i)\binom{k+1}{i} - (k-1-i)\binom{k}{i} - (k-i)\binom{k}{i-1} \\
=& (k-i)\binom{k}{i} - (k-1-i)\binom{k}{i} \\
=& \binom{k}{i}.
\end{align*}
Now, consider the case when $k$ is even and $i = \frac{k}{2}$. We have
\begin{align*}
[t^{\frac{k}{2}}] \left(Y_{\U_{k,k+1}}(t) - (1+t)\cdot Y_{\U_{k-1,k}}(t)\right) 
=& [t^{\frac{k}{2}}] Y_{\U_{k,k+1}}(t) - 2 \times [t^{\frac{k}{2}-1}] Y_{\U_{k-1,k}}(t) \\
=& (k - \frac{k}{2}) \binom{k+1}{\frac{k}{2}} - 2 \left(k - 1 - (\frac{k}{2} - 1)\right) \binom{k}{\frac{k}{2} - 1} \\
=& \binom{k}{\frac{k}{2}} \left( 1 - \frac{1}{1 + \frac{k}{2}} \right) \\
=& \binom{k}{\frac{k}{2}} - C_{\frac{k}{2}}.
\end{align*}
Thus, we complete the proof.
\qed

As an immediate consequence of Theorem \ref{sparsing_paving_matroid-inverse-Z}, we obtain the following result, which provides a more precise expression for the coefficients of the inverse $Z$-polynomials of sparse paving matroids.

\begin{corollary}\label{corollary-spas-paving-inver-Z-1}
 Let $\M$ be a sparse paving matroid of rank $k$ and cardinality $n$. Assume that $\M$ has exactly $\lambda$ circuit-hyperplanes. Then, the coefficients of $Y_{\M}(t)$ are given by
 \begin{align}\label{corollary-spas-paving-inver-Z}
 [t^i] Y_{\M}(t) =
 \begin{cases}
 \binom{n}{k} \binom{k}{i} \left( \frac{k-i}{n-i} - \lambda^* \right) & \text{for}~ 0 \leq i \leq \frac{k-1}{2}, \vspace{4mm} \\
 \binom{n}{k} \binom{k}{\frac{k}{2}} \cdot \frac{k}{2+k} \cdot \left( \frac{2+k}{2n-k} - \lambda^* \right) & \text{for}~ i = \frac{k}{2},
 \end{cases}
 \end{align}
 where $\lambda^* := \frac{\lambda}{\binom{n}{k}}$.
\end{corollary}

\proof
From Theorem \ref{sparsing_paving_matroid-inverse-Z}, for any $0 \leq i \leq \frac{k-1}{2}$, the coefficient of $t^i$ in $Y_{\M}(t)$ is
 \begin{align}
 [t^i] Y_{\U_{k,n}}(t) - \lambda \binom{k}{i} &= \binom{n}{i} \binom{n-i-1}{n-k} - \lambda \binom{k}{i} \label{another-form3} \\
 &= \frac{n!}{(n-k)!(k-i)!i!} \left( \frac{k-i}{n-i} - \frac{\lambda k!(n-k)!}{n!} \right) \nonumber \\
 &= \binom{n}{k} \binom{k}{i} \left( \frac{k-i}{n-i} - \frac{\lambda}{\binom{n}{k}} \right) \nonumber \\
 &= \binom{n}{k} \binom{k}{i} \left( \frac{k-i}{n-i} - \lambda^* \right). \label{SPP_IN_Z_CO_1}
 \end{align}

When $k$ is even and $i = \frac{k}{2}$, from Theorem \ref{sparsing_paving_matroid-inverse-Z}, the coefficient of $t^{\frac{k}{2}}$ is
\begin{align*}
 & \binom{n}{\frac{k}{2}} \binom{n - \frac{k}{2} - 1}{n-k} - \lambda \left[ \binom{k}{ \frac{k}{2}} - C_{\frac{k}{2}} \right]\\
 =& \binom{n}{\frac{k}{2}} \binom{n - \frac{k}{2} - 1}{n-k} - \lambda \cdot \frac{\frac{k}{2}}{1 + \frac{k}{2}} \binom{k}{ \frac{k}{2}} \\
 = &\frac{n!}{(n-k)! (\frac{k}{2})!(\frac{k}{2})! } 
 \left( \frac{\frac{k}{2}}{n - \frac{k}{2}} - \lambda \cdot \frac{k!(n-k)!}{n!} \cdot \frac{\frac{k}{2}}{1 + \frac{k}{2}} \right) \\
 = &\binom{n}{k} \binom{k}{\frac{k}{2}} \left( \frac{\frac{k}{2}}{n - \frac{k}{2}} - \frac{\lambda}{\binom{n}{k}} \cdot \frac{\frac{k}{2}}{1 + \frac{k}{2}} \right) \\
 =& \binom{n}{k} \binom{k}{\frac{k}{2}} \cdot \frac{\frac{k}{2}}{1 + \frac{k}{2}} \left( \frac{1 + \frac{k}{2}}{n - \frac{k}{2}} - \lambda^* \right).
\end{align*}
Thus, we complete the proof.
\qed

\begin{remark}\label{remark-spas-paving-inver-Z-1}
\textnormal{
When $k$ is an even number and $i = \frac{k}{2}$, we have the following expression for the coefficient of $t^{\frac{k}{2}}$ in $Y_{\M}(t)$, i.e.,
\begin{align*}
[t^{\frac{k}{2}}] Y_{\M}(t) &= \binom{n}{k} \binom{n - \frac{k}{2} - 1}{n - k} - \lambda \left[ \binom{k}{\frac{k}{2}} - C_{\frac{k}{2}} \right] \\
&= \binom{n}{k} \binom{n - \frac{k}{2} - 1}{n - k} - \lambda \binom{k}{\frac{k}{2}} + \lambda C_{\frac{k}{2}} \\
&= \binom{n}{k} \binom{k}{\frac{k}{2}} \left( \frac{\frac{k}{2}}{n - \frac{k}{2}} - \lambda^* \right) + \lambda C_{\frac{k}{2}}.
\end{align*}
Note that the first term, $\binom{n}{k} \binom{k}{\frac{k}{2}} \left( \frac{\frac{k}{2}}{n - \frac{k}{2}} - \lambda^* \right)$, is identical to the corresponding term in the first equality of \eqref{corollary-spas-paving-inver-Z} when $i = \frac{k}{2}$.
}
\end{remark}

\section{Unimodality and Log-Concavity}\label{4}

This section is dedicated to proving Theorem \ref{hjkl-1}, which demonstrates that the coefficients of the inverse $Z$-polynomials of sparse paving matroids are both unimodal and log-concave.

We have already established in Proposition \ref{inv_z_non_negative} that the coefficients of the inverse $Z$-polynomial $Y_{\M}(t)$ for any matroid $\M$ are non-negative. Based on this result, we now present a new proof of the non-negativity of the coefficients of the inverse $Z$-polynomials of sparse paving matroids, derived through direct calculation.

\proof[Second proof of Proposition \ref{inv_z_non_negative} in the case of sparse paving matroids]
From Corollary \ref{corollary-spas-paving-inver-Z-1} and Remark \ref{remark-spas-paving-inver-Z-1}, it suffices to show that
$$
\binom{n}{k}\binom{k}{i}\left(\frac{k-i}{n-i}-\lambda^*\right) \geq 0$$ 
for $0 \leq i \leq \frac{k}{2}.$
Note that for any two positive integers $2 \leq a \leq b$, we have 
$\frac{a}{b} \geq \frac{a-1}{b-1}.$
Thus, it follows that
\begin{align}\label{nongeative-sprsepaving-1}
\frac{k-i}{n-i} \geq \frac{k-i-1}{n-i-1} \geq \cdots \geq \frac{1}{n-k+1}.
\end{align}
On the other hand, from \cite[Corollary 4.13]{ferroni2022valuative}, we have
$\lambda \leq \binom{n}{k}\min \left\{\frac{1}{k+1}, \frac{1}{n-k+1}\right\},$
which implies
\begin{align}\label{nongeative-sprsepaving-2}
\lambda^* = \frac{\lambda}{\binom{n}{k}} \leq \min \left\{\frac{1}{k+1}, \frac{1}{n-k+1}\right\} \leq \frac{1}{n-k+1}.
\end{align}
Combining \eqref{nongeative-sprsepaving-1} and \eqref{nongeative-sprsepaving-2}, we obtain
\begin{align}\label{nongeative-sprsepaving-3}
\frac{k-i}{n-i} - \lambda^* \geq 0
\end{align}
for each $0 \leq i \leq \frac{k}{2}$.
Hence, the coefficients of the inverse $Z$-polynomials of sparse
paving matroids are non-negative.
\qed

Proudfoot, Xu, and Young \cite{PXY_2018ELC} proved that 
the coefficients of the $Z$-polynomial of a realizable matroid can be interpreted as the intersection cohomology Betti numbers of a projective variety. Later, in the work of Braden, Huh, Matherne, Proudfoot, and Wang \cite{BHM_ARX}, a more precise description of the intersection cohomology was obtained for any matroid $\M$. This naturally leads to the following question.

\textbf{Problem 1.} 
Do the coefficients of the inverse $Z$-polynomial $Y_{\M}(t)$ have a geometric interpretation?

Now, we proceed to prove Theorem \ref{hjkl-1},
which demonstrates that Conjecture \ref{conj-uni} and \ref{conj-log} hold for sparse paving matroids.

\proof
(1) Let us first prove that the coefficients of the inverse $Z$-polynomials of sparse paving matroids are unimodal. For $k=1$ or $k=2$, the case is trivial. We assume $k \geq 3$.
From Theorem \ref{sparsing_paving_matroid-inverse-Z}, the polynomial $Y_{\M}(t)$ is of degree $k$ and symmetric with center of symmetry $\lfloor k/2 \rfloor$. For $0 \leq i \leq k$, let 
$a_{i}:=[t^i]Y_{\M}(t)$
denote the coefficient of $t^i$ in $Y_{\M}(t)$. 

In order to prove the unimodality of $\{a_i\}_{i=0}^k$, we first show that $a_i \leq a_{i+1}$
for $0 \leq i \leq \frac{k-3}{2}$. 
To this end, we use the expression for $a_i$ from Corollary \ref{corollary-spas-paving-inver-Z-1}, i.e.,
\begin{align*}
a_i = [t^i] Y_{\M}(t) = \binom{n}{k} \binom{k}{i} \left( \frac{k-i}{n-i} - \lambda^* \right)
\end{align*}
for $0 \leq i \leq \frac{k-1}{2}$. Thus, for $0 \leq i \leq \frac{k-3}{2}$, it suffices to show
$$ \binom{n}{k} \binom{k}{i} \left( \frac{k-i}{n-i} - \lambda^* \right) \leq \binom{n}{k} \binom{k}{i+1} \left( \frac{k-i-1}{n-i-1} - \lambda^* \right), $$
which simplifies to
$$ \left[\binom{k}{i+1} - \binom{k}{i}\right] \lambda^* \leq \frac{k-i-1}{n-i-1} \binom{k}{i+1} - \frac{k-i}{n-i} \binom{k}{i}. $$
Note that
\begin{align*}
\binom{k}{i+1} - \binom{k}{i} &= \frac{k-2i-1}{k-i} \binom{k}{i+1} > 0,
\end{align*}
and
\begin{align*}
\frac{k-i-1}{n-i-1} \binom{k}{i+1} - \frac{k-i}{n-i} \binom{k}{i} &= \left( \frac{k-i-1}{n-i-1} - \frac{i+1}{n-i} \right) \binom{k}{i+1} \\
&= \frac{(n-i)(k-2i-2) + (i+1)}{(n-i)(n-i-1)} \binom{k}{i+1} > 0.
\end{align*}
Thus, it suffices to prove that for $0 \leq i \leq \frac{k-3}{2}$,
\begin{align*}
\lambda^* \leq \frac{(n-i)(k-2i-2)+(i+1)}{(n-i-1)(n-i)} \cdot \frac{k-i}{k-2i-1}.
\end{align*}
Recall that $\lambda^* \leq \min \left\{ \frac{1}{k+1}, \frac{1}{n-k+1} \right\} \leq \frac{1}{n-k+1}$, we now need to prove
$$ \frac{(n-i)(k-2i-2)+(i+1)}{(n-i)(n-i-1)} \cdot \frac{k-i}{k-2i-1} \geq \frac{1}{n-k+1}, $$ 
which is equivalent to
\begin{align*}
(n-k+1)(k-i) \left[ 2i^2 - (2n+k-3)i + n(k-2) + 1 \right] - (n-i)(n-i-1)(k-2i-1) \geq 0.
\end{align*}
This inequality simplifies to proving that for $0 \leq i \leq \frac{k-3}{2}$,
\begin{align*}
(k-n) \left[ 2i^3 - (2n+3k-5)i^2 + (3 - 4k + k^2 - 4n + 3kn)i + (1-k-n+3kn-k^2n) \right] \geq 0.
\end{align*}

Let 
$$ f(x) := 2x^3 - (2n+3k-5)x^2 + (3 - 4k + k^2 - 4n + 3kn)x + (1-k-n+3kn-k^2n) $$
be the continuous function on $[0, \frac{k-3}{2}]$. Then we compute its derivatives and obtain
\begin{align*}
f'(x) &= 6x^2 - 2(2n+3k-5)x + (3 - 4k + k^2 - 4n + 3kn), \\
f''(x) &= 12x - 2(2n+3k-5).
\end{align*}
Since $f''(x)$ increases monotonically on $[0, \frac{k-3}{2}]$ and $f''(\frac{k-3}{2}) = -2(2n+4) < 0$, it follows that $f''(x) \leq 0$ on $[0, \frac{k-3}{2}]$. Hence, $f'(x)$ decreases monotonically on $[0, \frac{k-3}{2}]$. From
$$ f'(\frac{k-3}{2}) = k(n - \frac{k}{2}) + 2n + k + \frac{3}{2} > 0 \quad \text{when} \quad k \geq 3, $$
we obtain that $f(x)$ increases monotonically on $[0, \frac{k-3}{2}]$. Combining this with
$$ f(\frac{k-3}{2}) = -\frac{1}{2}(kn + 2k - n - 2) = -\frac{1}{2}(n+2)(k-1) < 0, $$
we conclude that
$f(x) \leq 0 $ on $[0, \frac{k-3}{2}].$
Thus, we have
$$ 2i^3 - (2n+3k-5)i^2 + (3 - 4k + k^2 - 4n + 3kn)i + (1-k-n+3kn-k^2n) \leq 0 $$
for $0 \leq i \leq \frac{k-3}{2}$, which implies
$a_i \leq a_{i+1}.$

When $k$ is even and $i = \frac{k}{2} - 1$, we need to show
$a_{\frac{k}{2}-1} \leq a_{\frac{k}{2}}.$
From Corollary \ref{corollary-spas-paving-inver-Z-1}, we have
\begin{align*}
a_{\frac{k}{2}-1} &= \binom{n}{k} \binom{k}{\frac{k}{2}-1} \left( \frac{k+2}{2n-k+2} - \lambda^* \right), \\
a_{\frac{k}{2}} &= \binom{n}{k} \binom{k}{\frac{k}{2}} \cdot \frac{k}{k+2} \cdot \left( \frac{k+2}{2n-k} - \lambda^* \right).
\end{align*}
Thus, we need to show
$$ \binom{k}{\frac{k}{2}-1} \left( \frac{k+2}{2n-k+2} - \lambda^* \right) \leq \binom{k}{\frac{k}{2}} \cdot \frac{k}{k+2} \cdot \left( \frac{k+2}{2n-k} - \lambda^* \right). $$
Note that
$$ -\binom{k}{\frac{k}{2}-1} + \frac{k}{k+2} \binom{k}{\frac{k}{2}} = 0 $$
and
$$\frac{k}{k+2}\cdot \frac{k+2}{2n-k}\binom{k}{\frac{k}{2}}-\frac{k+2}{2n-k+2}\binom{k}{\frac{k}{2}-1}
=\frac{2k}{(2n-k+2)(2n-k)}\binom{k}{\frac{k}{2}}.$$
Hence, it suffices to show
$$ \frac{2k}{(2n-k+2)(2n-k)} \binom{k}{\frac{k}{2}} \geq 0, $$
which is obviously true. Thus, the proof is complete.

(2)
Now we proceed to prove the log-concavity of $\{a_i\}_{i=0}^k$. Let us first prove that for $1 \leq i < \frac{k}{2} - 1$,
\begin{align}\label{log_SPP_IN_Z_CO_11}
a_i^2 \geq a_{i-1}a_{i+1}.
\end{align}
To this end, we still use the expression for $a_i$ given by \eqref{corollary-spas-paving-inver-Z} in Corollary \ref{corollary-spas-paving-inver-Z-1} for $0 \leq i < \frac{k}{2}$. 
Thus, \eqref{log_SPP_IN_Z_CO_11} becomes
\begin{align*}
\binom{k}{i}^2 \left( \frac{k-i}{n-i} - \lambda^* \right)^2 \geq \binom{k}{i-1} \binom{k}{i+1} \left( \frac{k-i+1}{n-i+1} - \lambda^* \right) \left( \frac{k-i-1}{n-i-1} - \lambda^* \right),
\end{align*}
for $1 \leq i < \frac{k}{2} - 1$. 
Note that
$
\binom{k}{i}^2 > \binom{k}{i-1} \binom{k}{i+1}
$
for $1 \leq i \leq k-1$. It suffices to show that for any $1 \leq i < \frac{k}{2} - 1$,
\begin{align*}
\left( \frac{k-i}{n-i} - \lambda^* \right)^2 \geq \left( \frac{k-i+1}{n-i+1} - \lambda^* \right) \left( \frac{k-i-1}{n-i-1} - \lambda^* \right).
\end{align*}
This is equivalent to proving
\begin{align}\label{log_SPP_IN_Z_CO_2}
\left( \frac{k-i+1}{n-i+1} + \frac{k-i-1}{n-i-1} - \frac{2(k-i)}{n-i} \right) \lambda^* \geq \frac{(k-i+1)(k-i-1)}{(n-i+1)(n-i-1)} - \frac{(k-i)^2}{(n-i)^2}.
\end{align}
Note that
\begin{align*}
\frac{k-i+1}{n-i+1} + \frac{k-i-1}{n-i-1} - \frac{2(k-i)}{n-i} 
&= 2 \left[ \frac{(n-i)(k-i)-1}{(n-i)^2-1} - \frac{(n-i)(k-i)}{(n-i)^2} \right] \\ \vspace{1mm}
& = -\frac{2(n-k)(n-i)}{(n-i)^4 - (n-i)^2},
\end{align*}
and
\begin{align*}
\frac{(k-i+1)(k-i-1)}{(n-i+1)(n-i-1)} - \frac{(k-i)^2}{(n-i)^2} 
= \frac{(k-i)^2 - 1}{(n-i)^2 - 1} - \frac{(k-i)^2}{(n-i)^2}
= -\frac{(n-k)(n+k-2i)}{(n-i)^4 - (n-i)^2}.
\end{align*}
Thus, we need to show that for any $1 \leq i < \frac{k}{2} - 1$,
$$
\lambda^* \leq \frac{n+k-2i}{2(n-i)}.
$$
In fact, combining with
$2 \leq n+k-2i \leq 2(n-i)$ and \eqref{nongeative-sprsepaving-3}, we have
$$
\frac{n+k-2i}{2(n-i)} \geq \frac{k-i}{n-i} \geq \lambda^*,
$$
as desired.

When $k$ is odd and $i = \frac{k-1}{2}$, we need to prove
\begin{align}\label{log_SPP_IN_Z_CO_33}
a_{\frac{k-1}{2}}^2 \geq a_{\frac{k-3}{2}} a_{\frac{k+1}{2}}.
\end{align}
From Theorem \ref{sparsing_paving_matroid-inverse-Z} and Theorem \ref{hjkl-1}, we have
$$
a_{\frac{k-1}{2}} = a_{\frac{k+1}{2}}, \quad a_{\frac{k-1}{2}} \geq a_{\frac{k-3}{2}}.
$$
Thus, \eqref{log_SPP_IN_Z_CO_33} is an obvious fact.

When $k$ is even and $i = \frac{k}{2}-1$, we need to prove
\begin{align}\label{k_even-log}
a_{\frac{k}{2}-1}^2 \geq a_{\frac{k}{2}-2} a_{\frac{k}{2}}.
\end{align}
From Theorem \ref{sparsing_paving_matroid-inverse-Z}, we have
\begin{align*}
a_{\frac{k}{2}-2} &= \binom{n}{k} \binom{k}{\frac{k}{2}-2} \left( \frac{k+4}{2n-k+4} - \lambda^* \right), \\
a_{\frac{k}{2}-1} &= \binom{n}{k} \binom{k}{\frac{k}{2}-1} \left( \frac{k+2}{2n-k+2} - \lambda^* \right), \\
a_{\frac{k}{2}} &= \binom{n}{k} \binom{k}{\frac{k}{2}} \cdot \frac{k}{k+2} \left( \frac{k+2}{2n-k} - \lambda^* \right).
\end{align*}
Thus, \eqref{k_even-log} becomes
\begin{align*}
\binom{k}{\frac{k}{2}-1}^2 \left( \frac{k+2}{2n-k+2} - \lambda^* \right)^2 \geq \frac{k}{k+2} \binom{k}{\frac{k}{2}-2} \binom{k}{\frac{k}{2}} \left( \frac{k+4}{2n-k+4} - \lambda^* \right) \left( \frac{k+2}{2n-k} - \lambda^* \right),
\end{align*}
or equivalently,
\begin{align} \label{final-log-sp}
\left( \frac{k+2}{2n-k+2} - \lambda^* \right)^2 \geq \frac{k-2}{k+4} \left( \frac{k+4}{2n-k+4} - \lambda^* \right) \left( \frac{k+2}{2n-k} - \lambda^* \right).
\end{align}
To complete the proof, we verify that \eqref{final-log-sp} holds for any positive integers $1 \leq k \leq n$. While this could be confirmed through a tedious computation, we prefer to provide a computer-aided proof as follows.

\begin{mma}
\In \left( \frac{k+2}{2n-k+2} - \lambda^* \right)^2 \geq \frac{k-2}{k+4} \left( \frac{k+4}{2n-k+4} - \lambda^* \right) \left( \frac{k+2}{2n-k} - \lambda^* \right); \\ 
\In 0\leq k\leq n \&\& n\geq 1\&\&
0\leq\lambda^*\leq |Min|\left[\frac {1} {k + 1}, \frac {1} {-k + n + 1} \right]; \\
\In 
|Resolve| \left[|ForAll| [ {
\left \{k, n, \lambda^*\right\} \ }, \%,\%\%]\right] \\
\Out \text{True} \\
\end{mma}

When $k$ is even and $i=\frac{k}{2}$, we need to prove
\begin{align*}
a_{\frac{k}{2}}^2 \geq a_{\frac{k}{2}-1}a_{\frac{k}{2}+1}.
\end{align*} 
Form Theorem \ref{sparsing_paving_matroid-inverse-Z} and Theorem \ref{hjkl-1}, we have 
$$ a_{\frac{k}{2}-1}= a_{\frac{k}{2}+1}~~~~\text{and}~~~~ a_{\frac{k}{2}}\geq a_{\frac{k}{2}-1}. $$
Thus, the proof is complete.
\qed

A property known as $\gamma$-positivity, which implies symmetry and unimodality, has recently attracted attention in topological, algebraic, and enumerative combinatorics. We recall that a polynomial $f(t) = \sum_i a_i t^i \in \mathbb{R}[t]$ is said to be $\gamma$-positive if it can be expressed as
$$
f(t) = \sum_{i=0}^{\lfloor \frac{n}{2} \rfloor} \gamma_i t^i (1+t)^{n-2i},
$$
where $n \in \mathbb{N}$ and $\gamma_0, \gamma_1, \ldots, \gamma_{\lfloor \frac{n}{2} \rfloor}$ are non-negative real numbers. Ferroni and Schröter \cite{ferroni2022valuative} proved that the $Z$-polynomials of matroids are $\gamma$-positive. Since $Y_{\M}(t)$ is a palindromic polynomial, it is natural to investigate its $\gamma$-positive property.
However, we argue that the inverse $Z$-polynomial $Y_{\U_{k,n}}(t)$ of a matroid is not necessarily $\gamma$-positive. As an example, consider the uniform matroid $\U_{4,5}$. According to Theorem \ref{main-theorem-uniform}, the inverse $Z$-polynomial of $\U_{4,5}$ is given by
$Y_{\U_{4,5}}(t) = 4t^4 + 15t^3 + 20t^2 + 15t + 4.$
Upon calculation, we find that the associated $\gamma$-polynomial is
$\gamma(Y_{\U_{4,5}}(t), t) = -2t^2 - t + 4.$
Since the coefficients of this $\gamma$-polynomial are not all non-negative, it follows that $Y_{\U_{4,5}}(t)$ is not $\gamma$-positive.

\vspace{4mm}

\noindent{\bf Acknowledgements.} 

The first author is supported by the National Science Foundation of China (No.11801447) and Guangdong Basic and Applied Basic Research Foundation (No.2024A1515011276).
The third author is supported by the National Science Foundation of China (No.11901431).

\bibliographystyle{plain}

\end{document}